\documentclass[bezier]{article}
\usepackage{graphicx}
\usepackage{amsmath,amssymb,amsfonts, euscript}
\usepackage{multirow}
\usepackage{mathtools}
\newtheorem{theorem}{Theorem}[section]

\newtheorem{example}[theorem]{Example}

\newtheorem{algorithm}[theorem]{Algorithm}

\title{\bf\Large Minimal residual Hermitian and skew-Hermitian splitting iteration method for the continuous Sylvester equation}

\author{Zeinab Bahramizadeh$^1$, Mojtaba Nazari$^1$, Mohammad Khorsand Zak$^2$\footnote{Corresponding author, {\it email}: mo.khorsand@mail.um.ac.ir} and Zahra Yarahmadi$^1$
\\{\footnotesize{\it $^1$Department of Mathematics, Khoramabad Branch, IAU, Khoramabad, Iran}}
\\{\footnotesize{\it $^2$Department of Mathematics, Aligudarz Branch, IAU, Aligudarz, Iran}}
}

\date{}

\begin{document}
\maketitle

\begin{abstract}
By applying the minimal residual technique to the Hermitian and skew-Hermitian (HSS) iteration scheme, we introduce a non-stationary iteration method named minimal residual Hermitian and skew-Hermitian (MRHSS) iteration method to solve the continuous Sylvester equation. Numerical results verify the effectiveness and robustness of the MRHSS iteration method versus the HSS method for the continuous Sylvester equation. Moreover, by numerical computation, we show that the MRHSS splitting can be used as a splitting preconditioner and induce accurate, robust and effective preconditioned Krylov subspace iteration methods for solving the continuous Sylvester equation.

\end{abstract}

\section{Introduction}
In many problems in scientific computing we encounter with matrix equations. Matrix equations are one of the most interesting and intensively studied classes of mathematical problems and play vital roles in applications, and many researchers have studied matrix equations and their applications, see \cite{Datta, Dehghan, Dehghan1, Hajarian, Ke, Khorsand1, Simoncini, Tohidi} and their references.
 Nowadays, the continuous Sylvester equation is possibly the most famous and the most broadly employed linear matrix equation, and is given as
\begin{equation}\label{1}
AX+XB=C,
\end{equation}
where $A\in\mathbb{C}^{n\times n}$, $B\in\mathbb{C}^{m\times m}$ and $C\in\mathbb{C}^{n\times m}$ are defined matrices and $X\in\mathbb{C}^{n\times m}$ is an unknown matrix. A Lyapunov equation is a special case with $m=n, B=A^T$, and $C=C^T$. Here and in the sequel, $W^T$ is used to denote the transpose of the matrix $W\in\mathbb{R}^{n\times n}$. Equation \eqref{1} has a unique solution if and only if $A$ and $-B$ have no common eigenvalues, which will be assumed throughout this paper.

Many results have been obtained about the Sylvester equation and it appears frequently in many areas of applied mathematics and plays vital roles in a number of applications such as control theory \cite{Datta}, model reduction \cite{Benner} and image processing \cite{Bouhamidi}, see \cite{Bai, Bartels, Dehghan, Evans, El, Golub, Hu, Khorsand2, Khorsand3, Salkuyeh} and their references for more details.

In general, the dimensions of $A$ and $B$ may be orders of magnitude different, and this fact is key in selecting the most appropriate numerical solution strategy \cite{Simoncini}.
For solving general Sylvester equations of small size we use some methods which classified such as direct methods. Some of these direct methods are the Bartels-Stewart \cite{Bartels} and the Hessenberg-Schur \cite{Golub} methods which consist of transforming coefficient matrices $A$ and $B$ into triangular or Hessenberg form by an orthogonal similarity transformation and then solving the resulting system directly by a back-substitution process.
When the coefficient matrices $A$ and $B$ are large and sparse, iterative methods are often the methods of choice for solving the Sylvester equation \eqref{1} efficiently and accurately. Many iterative methods were developed for solving matrix equations, such as the alternating direction implicit (ADI) method \cite{Benner}, the Krylov subspace based algorithms \cite{Hu, Salkuyeh, El}, the Hermitian and skew-Hermitian splitting (HSS) method, and the inexact variant of HSS (IHSS) iteration method \cite{Bai2}, The nested splitting conjugate gradient (NSCG) method \cite{Khorsand2} and the nested splitting CGNR (NS-CGNR) method \cite{Khorsand3}.

When both coefficient matrices are (non-Hermitian) positive semi-definite, and at least one of them is positive definite, the Hermitian and skew-Hermitian splitting (HSS) method \cite{Bai} and the nested splitting conjugate gradient (NSCG) method \cite{Khorsand2} are often the methods of choice for efficiently and accurately solving the Sylvester equation \eqref{1}.

In order to study the numerical methods, we often rewrite the continuous Sylvester equation \eqref{1} as a mathematically equivalent linear system of equations such as follows:
\begin{equation}\label{2}
\mathcal{A}x=c,
\end{equation}
where the matrix $\mathcal{A}$ is of dimension $nm\times nm$ and is given by
\begin{equation}\label{3}
\mathcal{A}=I_m\otimes A+B^T\otimes I_n,
\end{equation}
where $\otimes$ denotes the Kronecker product $(A\otimes B=[a_{ij}B])$ and
$$
\begin{array}{l}
c=vec(C)=(c_{11},c_{21},\cdots,c_{n1},c_{12},c_{22},\cdots,c_{n2},\cdots,c_{nm})^T,\\
x=vec(X)=(x_{11},x_{21},\cdots,x_{n1},x_{12},x_{22},\cdots,x_{n2},\cdots,x_{nm})^T.
\end{array}
$$

Of course, this is a numerically poor way to determine the solution $X$ of the Sylvester equation \eqref{1}, as the linear system of equations \eqref{2} is costly to solve and can be ill-conditioned.

Motivated by \cite{Yang1, Yang2}, we apply the minimal residual technique to the Hermitian and skew-Hermitian iteration scheme and introduce a non-stationary iteration method named minimal residual Hermitian and skew-Hermitian (MRHSS) iteration method to solve the continuous Sylvester equation.

In the remainder of this paper, we use $||M||_2$, $||M||_F$ and $I_n$ to denote the spectral norm, the Frobenius norm of a matrix $M\in\mathbb{C}^{n\times n}$, and the identity matrix with dimension $n$, respectively. Note that $||.||_2$ is also used to represent the 2-norm of a vector.
Furthermore, we have the following equivalent relationships between the Frobenius norm of a matrix $R$ and the 2-norm of a vector $r=vec(R)$:
$$||r||_2=\sqrt{\sum_{i=1}^{mn}|r_i|^2}\Leftrightarrow||R||_F=\sqrt{\sum_{i=1}^m\sum_{j=1}^n|R_{ij}|^2}.$$

The reminder of this paper is organized as follows. Section \ref{main} presents the minimal residual Hermitian and skew-Hermitian splitting (MRHSS) method for the continuous Sylvester equation. Section \ref{Numer} is devoted to numerical experiments. Finally, we present our conclusions in Section \ref{conc}.

\section{Main results}\label{main}
For the linear system of equations \eqref{2}, we consider the Hermitian and skew-Hermitian splitting $\mathcal{A}=\mathcal{H}+\mathcal{S}$, where
\begin{equation}\label{5}
 \mathcal{H}=\frac{\mathcal{A}+\mathcal{A}^T}{2}, ~~~\mathcal{S}=\frac{\mathcal{A}-\mathcal{A}^T}{2},
\end{equation}
are the Hermitian and skew-Hermitian parts of matrix $\mathcal{A}$, respectively. Then, the iteration scheme of the MRHSS iteration method \cite{Yang1, Yang2} for system of linear equations \eqref{2} is
\begin{equation}\label{e1}
  \left\{
  \begin{array}{rcl}
    x^{(k+\frac{1}{2})} & = & x^{(k)}+\beta_k\delta^{(k)} \\
    x^{(k+1)} & = & x^{(k+\frac{1}{2})}+\gamma_k\delta^{(k+\frac{1}{2})},
  \end{array}
  \right.
\end{equation}
where, $\delta^{(k)}=(\hat{\alpha} I+\mathcal{H})^{-1}r^{(k)}$, $\delta^{(k+\frac{1}{2})}=(\hat{\alpha} I+\mathcal{S})^{-1}r^{(k+\frac{1}{2})}$,
$r^{(k)}=c-\mathcal{A}x^{(k)}$ and $r^{(k+\frac{1}{2})}=c-\mathcal{A}x^{(k+\frac{1}{2})}$. Let $M_1=\mathcal{A}(\hat{\alpha} I+\mathcal{H})^{-1}$ and $M_2=\mathcal{A}(\hat{\alpha} I+\mathcal{S})^{-1}$. The residual form of iteration scheme \eqref{e1} can be written as
\begin{equation}\label{e2}
  \left\{
  \begin{array}{rcl}
    r^{(k+\frac{1}{2})} & = & r^{(k)}-\beta_kM_1r^{(k)} \\
    r^{(k+1)} & = & r^{(k+\frac{1}{2})}-\gamma_kM_2r^{(k+\frac{1}{2})}.
  \end{array}
  \right.
\end{equation}
Denote $M=(\hat{\alpha}I+\mathcal{H})^{-1}$. Then, an inner product can be defined as
\begin{equation}\label{e3}
  (x,y)_M=(Mx,My), \hspace{1cm}\forall x,y\in \mathbb{C}^{nm},
\end{equation}
where $(\cdot,\cdot)$ denotes the $l^2$ inner product of two vectors. Thus, for $x\in\mathbb{C}^{nm}$ and $X\in\mathbb{C}^{nm\times nm}$, the induced vector and the induced matrix norms can be defined as $||x||_M=||Mx||_2$ and $||X||_M=||MXM^{-1}||_2$, respectively. Now, the parameter $\beta_k$ is determined by the 2-norm of the residual, and we have
\begin{equation}\label{e4}
  \beta_k=\frac{(r^{(k)},M_1r^{(k)})}{||M_1r^{(k)}||_2^2}.
\end{equation}
However, the parameter $\gamma_k$ will be determined by minimizing the M-norm of the residual rather than the 2-norm, see \cite{Yang1}. Therefore, we have
\begin{equation}\label{e5}
  \gamma_k=\frac{(Mr^{(k+\frac{1}{2})},MM_2r^{(k+\frac{1}{2})})}{||MM_2r^{(k+\frac{1}{2})}||_2^2}.
\end{equation}
According the following theorem, the iteration scheme \eqref{e1} is an unconditionally convergent MRHSS iteration method \cite{Yang1}.
\begin{theorem}\label{thm1}
Let $\mathcal{A}$ be a non-Hermitian positive definite matrix. Then, the MRHSS iteration method used for solving the system of linear equations \eqref{2} is unconditionally convergent for any $\hat{\alpha}>0$ and any initial guess $x^{(0)}\in\mathbb{C}^{mn}$.
\end{theorem}
\textrm{Proof.} See \cite{Yang1}.

Let $H_A(\alpha)=\alpha I_n+H_A, S_A(\alpha)=\alpha I_n+S_A, H_{B^T}(\alpha)=\alpha I_m+H_{B^T}, S_{B^T}(\alpha)=\alpha I_m+S_{B^T}$ and $H_A,S_A,H_{B^T},S_{B^T}$ are the Hermitian and skew-Hermitian parts of $A$ and $B^T$, respectively.
For the Sylvester equation \eqref{1}, according to iterative scheme \eqref{e1}, we have the following iteration scheme
\begin{equation}\label{e6}
  \left\{
  \begin{array}{rcl}
    X^{(k+\frac{1}{2})} & = & X^{(k)}+\beta_k\Delta^{(k)} \\
    X^{(k+1)} & = & X^{(k+\frac{1}{2})}+\gamma_k\Delta^{(k+\frac{1}{2})},
  \end{array}
  \right.
\end{equation}
where, $\Delta^{(0)}$ obtain from the Sylvester equation
\begin{equation}\label{e7}
  H_A(\alpha)\Delta^{(0)}+\Delta^{(0)}H_B(\alpha)=R^{(0)},
\end{equation}
and $\Delta^{(k+\frac{1}{2})}$ obtain from the Sylvester equation
\begin{equation}\label{e8}
  S_A(\alpha)\Delta^{(k+\frac{1}{2})}+\Delta^{(k+\frac{1}{2})}S_B(\alpha)=R^{(k+\frac{1}{2})},
\end{equation}
with $R^{(0)}=C-AX^{(0)}-X^{(0)}B$ and $R^{(k+\frac{1}{2})}=C-AX^{(k+\frac{1}{2})}-X^{(k+\frac{1}{2})}B$.
We state how to update $\Delta^{(k+1)}$ a few later.

If the Sylvester equation \eqref{1} has a unique solution, then under the assumption $A$ and $B$ are positive semi-definite and at last one of them is positive definite, we can easily see that there is no common eigenvalue between the matrices $H_A$ and $-H_B$ (also for $S_A$ and $-S_B$), so the Sylvester equations \eqref{e7} and \eqref{e8} have unique solution for all given right hand side matrices.

From \eqref{3} and \eqref{5}, by using the Kronecker product's properties, we have
\begin{equation}\label{6}
\hat{\alpha}I+\mathcal{H}=I_m\otimes H_A(\alpha)+H_{B^T}(\alpha)\otimes I_n
\end{equation}
\begin{equation}\label{7}
\hat{\alpha}I+\mathcal{S}=I_m\otimes S_A(\alpha)+S_{B^T}(\alpha)\otimes I_n
\end{equation}
where $\alpha=\frac{\hat{\alpha}}{2}$. Form relations \eqref{e2}, we can obtain
\begin{equation}\label{e9}
\left\{
\begin{array}{rcl}
  R^{(k+\frac{1}{2})} & = & R^{(k)}-\beta W^{(k)} \\
  R^{(k+1)} & = & R^{(k+\frac{1}{2})}-\gamma W^{(k+\frac{1}{2})}
\end{array}
\right.
\end{equation}
where $W^{(k)}=A\Delta^{(k)}+\Delta^{(k)}B$ and $W^{(k+\frac{1}{2})}=A\Delta^{(k+\frac{1}{2})}+\Delta^{(k+\frac{1}{2})}B$. Moreover, similar to
\eqref{e4} and \eqref{e5}, we can obtain
\begin{equation}\label{e10}
 \beta=\frac{\langle R^{(k)},W^{(k)}\rangle_F}{\langle W^{(k)},W^{(k)}\rangle_F},
\end{equation}
and
\begin{equation}\label{e11}
 \gamma=\frac{\langle V^{(k+\frac{1}{2})},U^{(k+\frac{1}{2})}\rangle_F}{\langle U^{(k+\frac{1}{2})},U^{(k+\frac{1}{2})}\rangle_F},
\end{equation}
where, $V^{(k+\frac{1}{2})}$ obtain from the Sylvester equation
\[
H_A(\alpha)V^{(k+\frac{1}{2})}+V^{(k+\frac{1}{2})}H_B(\alpha)=R^{(k+\frac{1}{2})},
\]
and $U^{(k+\frac{1}{2})}$ obtain from the Sylvester equation
\[
H_A(\alpha)U^{(k+\frac{1}{2})}+U^{(k+\frac{1}{2})}H_B(\alpha)=W^{(k+\frac{1}{2})}
\]

On the surface, four systems of linear equations should be solved at each step of the MRHSS method for system of linear equations \eqref{2}. But it can be reduced to three. Denote $\zeta^{(k+\frac{1}{2})}=(\hat{\alpha}I+\mathcal{H})^{-1}r^{(k+\frac{1}{2})}$ and $v^{(k+\frac{1}{2})}=(\hat{\alpha}I+\mathcal{H})^{-1}\mathcal{A}\delta^{(k+\frac{1}{2})}$, the vector $\delta^{(k+1)}$ in Step $k+1$ can be calculated as follows
\[
\begin{array}{rcl}
\delta^{(k+1)}& = & (\hat{\alpha}I+\mathcal{H})^{-1}(c-\mathcal{A}x^{(k+1)}) \\
    & = & (\hat{\alpha}I+\mathcal{H})^{-1}(c-\mathcal{A}(x^{(k+\frac{1}{2})}+\gamma_k\delta^{(k+\frac{1}{2})})) \\
    & = & (\hat{\alpha}I+\mathcal{H})^{-1}(r^{(k+\frac{1}{2})}-\gamma_k\mathcal{A}\delta^{(k+\frac{1}{2})}) \\
    & = & \zeta^{(k+\frac{1}{2})}-\gamma_kv^{(k+\frac{1}{2})},
\end{array}
\]
where the $\zeta^{(k+\frac{1}{2})}$ and $v^{(k+\frac{1}{2})}$ have been calculated in Step $k$. Therefore, in \eqref{e6} we can update $\Delta^{(k+1)}$ as
\[
\Delta^{(k+1)}=V^{(k+\frac{1}{2})}-\gamma U^{(k+\frac{1}{2})}.
\]
In addition, we choose the value of parameter $\alpha$ as in \cite{Bai}.

Therefore, an implementation of the MRHSS method for the continuous Sylvester equation can be given by the following algorithm.

\begin{algorithm}\label{alg1}
\textbf{The MRHSS algorithm for the Sylvester equation}
\begin{itemize}
\item[1.]
Select an initial guess $X^{(0)}$, compute $R^{(0)}=C-AX^{(0)}-X^{(0)}B$
\item[2.]
Solve $H_A(\alpha)\Delta^{(0)}+\Delta^{(0)}H_B(\alpha)=R^{(0)}$
\item[3.]
For $k=0,1,2,\cdots,$ until convergence, Do:
\item[4.]
\hspace{0.7cm}$W^{(k)}=A\Delta^{(k)}+\Delta^{(k)}B$
\item[5.]
\hspace{0.7cm}$\beta=\frac{\langle R^{(k)},W^{(k)}\rangle_F}{\langle W^{(k)},W^{(k)}\rangle_F}$
\item[6.]
\hspace{0.7cm}$X^{(k+\frac{1}{2})}=X^{(k)}+\beta\Delta^{(k)}$
\item[7.]
\hspace{0.7cm}$R^{(k+\frac{1}{2})}=R^{(k)}-\beta W^{(k)}$
\item[8.]
\hspace{0.7cm}Solve $S_A(\alpha)\Delta^{(k+\frac{1}{2})}+\Delta^{(k+\frac{1}{2})}S_B(\alpha)=R^{(k+\frac{1}{2})}$
\item[9.]
\hspace{0.7cm}Solve $H_A(\alpha)V^{(k+\frac{1}{2})}+V^{(k+\frac{1}{2})}H_B(\alpha)=R^{(k+\frac{1}{2})}$
\item[10.]
\hspace{0.6cm}$W^{(k+\frac{1}{2})}=A\Delta^{(k+\frac{1}{2})}+\Delta^{(k+\frac{1}{2})}B$
\item[11.]
\hspace{0.6cm}Solve $H_A(\alpha)U^{(k+\frac{1}{2})}+U^{(k+\frac{1}{2})}H_B(\alpha)=W^{(k+\frac{1}{2})}$
\item[12.]
\hspace{0.6cm}$\gamma=\frac{\langle V^{(k+\frac{1}{2})},U^{(k+\frac{1}{2})}\rangle_F}{\langle U^{(k+\frac{1}{2})},U^{(k+\frac{1}{2})}\rangle_F}$
\item[13.]
\hspace{0.6cm}$X^{(k+1)}=X^{(k+\frac{1}{2})}+\gamma\Delta^{(k+\frac{1}{2})}$
\item[14.]
\hspace{0.6cm}$R^{(k+1)}=R^{(k+\frac{1}{2})}-\gamma W^{(k+\frac{1}{2})}$
\item[15.]
\hspace{0.6cm}$\Delta^{(k+1)}=V^{(k+\frac{1}{2})}-\gamma U^{(k+\frac{1}{2})}$
\item[16.]
End Do
\end{itemize}
\end{algorithm}

\begin{theorem}\label{thm2}
Suppose that the coefficient matrices $A$ and $B$ in the continuous Sylvester equation \eqref{1} are non-Hermitian positive semi-definite, and at least one of them is positive definite. Then the MRHSS iteration method \eqref{e6} for solving the Sylvester equation \eqref{1} is unconditionally convergent for any $\alpha>0$ and any initial guess $X^{(0)}\in\mathbb{C}^{n\times m}$.
\end{theorem}
\textrm{Proof.}
The continuous Sylvester equation \eqref{1} is mathematically equivalent to the linear system of equations \eqref{2}. Therefore, the proof is similar to that of Theorem 3.3 in \cite{Yang1} with only technical modifications.

\subsection{Using the MRHSS splitting as a preconditioner}
From the fact that any matrix splitting can naturally induce a splitting preconditioner for the Krylov subspace methods (see \cite{Bai2}) in section \ref{Numer}, by numerical computation, we show that the minimal residual Hermitian and skew-Hermitian splitting can be used as a splitting preconditioner and induce accurate, robust and effective preconditioned Krylov subspace iteration methods for solving the continuous Sylvester equation.

\section{Numerical results}\label{Numer}
In this section, we use a few numerical results to show the effectiveness of the MRHSS method by comparing its results with the HSS method. All numerical experiments were computed in double precision with a number of MATLAB codes. All iterations are started from the zero matrix for initial $X^{(0)}$ and terminated when the current iterate satisfies $$\frac{\|R^{(k)}\|_F}{\|R^{(0)}\|_F}\leq10^{-8},$$ where $R^{(k)}=C-AX^{(k)}-X^{(k)}B$ is the residual of the $k$th iterate. Also, we use the tolerance $\varepsilon=0.001$ for inner iterations in corresponding methods.
We report the results of the CPU time (CPU), the number of iteration steps (IT) and the norm of residual $\|R^{(k)}\|_F$ (res-norm) in the tables, and compare the HSS iterative method \cite{Bai} with the MRHSS iterative method for solving the continuous Sylvester equation \eqref{1}.

\begin{example}\label{Ex1}
 For this example, we use the coefficient matrices
\[
A=M+2rN+\frac{100}{(n+1)^2}I,\hspace{1cm}\text{and}\hspace{1cm}B=M+2rN+\frac{100}{(m+1)^2}I
\]
where $M={\rm tridiag}(-1,2,-1)$, $N={\rm tridiag}(0.5,0,-0.5)$ from suitable dimensions, and $r=0.01$ \cite{Bai, Khorsand1}.
\end{example}
This class of problems may arise in the preconditioned Krylov subspace iteration methods used for solving the systems of linear equations resulting from the finite difference or Sinc-Galerkin discretization of various differential equations and boundary value problems \cite{Bai}.

We apply the iteration methods to this problem with different dimensions $(n,m)$. The results are given in Tables \ref{Tab1a} and \ref{Tab1b}. From the results presented in the Tables \ref{Tab1a} and \ref{Tab1b}, we observe that the MRHSS method is more efficient than the HSS method in terms of CPU time. However, when the dimension increases, we observe that the HSS method is more efficient than the MRHSS method in terms of number of iterations (IT).
\begin{table}[h!]
\centering
\caption{The results for the example \ref{Ex1}}
\label{Tab1a}
\begin{tabular}{cccccccc}
\hline
\multirow{2}{*}{}
     & \multicolumn{3}{c}{HSS} &  & \multicolumn{3}{c}{MRHSS} \\
\cline{2-4} \cline{6-8}
 $(n,m)$ & CPU   & iteration& res-norm&  & CPU    & iteration& res-norm        \\
\hline
 $(8,8)$    & 0.04   & 14   & 2.3191e-6  &  & 0.02   & 7   & 2.3518e-6    \\
 $(16,16)$  & 0.05   & 26   & 1.2712e-6  &  & 0.03   & 16  & 1.3088e-6    \\
 $(32,32)$  & 0.16   & 48   & 1.3215e-6  &  & 0.12   & 37  & 1.1597e-6    \\
 $(64,64)$  & 1.02   & 89   & 1.5946e-6  &  & 0.91   & 85  & 1.6722e-6    \\
$(128,128)$ & 13.09  & 164  & 2.2369e-6  &  & 11.51  & 188 & 2.2271e-6    \\
$(256,256)$ & 85.04  & 298  & 3.2107e-6  &  & 75.06  & 404 & 3.2155e-6    \\
 \hline
\end{tabular}
\end{table}

\begin{table}[h!]
\centering
\caption{The results for the example \ref{Ex1}}
\label{Tab1b}
\begin{tabular}{cccccccc}
\hline
\multirow{2}{*}{}
     & \multicolumn{3}{c}{HSS} &  & \multicolumn{3}{c}{MRHSS} \\
\cline{2-4} \cline{6-8}
 $(n,m)$ & CPU   & iteration& res-norm&  & CPU    & iteration& res-norm        \\
\hline
 $(512,8)$  & 0.95   & 20   & 6.9889e-6  &  & 0.12   & 11  & 6.8967e-6    \\
 $(512,16)$ & 2.64   & 36   & 5.0093e-6  &  & 0.71   & 24  & 4.9071e-6    \\
 $(512,32)$ & 6.95   & 67   & 3.6776e-6  &  & 2.56   & 53  & 3.1928e-6    \\
 $(512,64)$ & 25.01  & 122  & 3.4599e-6  &  & 9.73   & 126 & 3.2791e-6    \\
$(512,128)$ & 90.23  & 218  & 3.4718e-6  &  & 39.60  & 272 & 3.5181e-6    \\
$(512,256)$ & 370.45 & 365  & 3.8891e-6  &  & 206.07 & 517 & 3.9374e-6    \\
 \hline
\end{tabular}
\end{table}

\begin{example}\label{Ex2}
We consider the continuous Sylvester equation \eqref{1} with the coefficient matrices
\[
\left\{
\begin{array}{l}
  A=\text{diag}(1, 2, \cdots, n)+rL^T, \\
  B=2^{-t}I_n+\text{diag}(1, 2, \cdots, n)+rL^T+2^{-t}L,
\end{array}
\right.
\]
with $L$ the strictly lower triangular matrix having ones in the lower triangle part \cite{Bai}. Here, $t$ is a problem parameter to be specified in actual computations.
\end{example}
The results of this problem are given in Table \ref{Tab2}. Here, we observe that the MRHSS method is more efficient in both terms of CPU time and number of iterations (IT) than the HSS method.
\begin{table}[h!]
\centering
\caption{The results for the example \ref{Ex2}}
\label{Tab2}
\begin{tabular}{cccccccc}
\hline
\multirow{2}{*}{}
     & \multicolumn{3}{c}{HSS} &  & \multicolumn{3}{c}{MRHSS} \\
\cline{2-4} \cline{6-8}
 $(n,m)$    & CPU    & IT   & res-norm   &  & CPU    & IT  & res-norm     \\
\hline
 $(8,8)$    & 0.04   & 19   & 6.9896e-6  &  & 0.02   & 11  & 7.3379e-6    \\
 $(16,16)$  & 0.07   & 24   & 1.7183e-5  &  & 0.03   & 16  & 2.2925e-5    \\
 $(32,32)$  & 0.14   & 31   & 8.1598e-5  &  & 0.07   & 22  & 9.4150e-5    \\
 $(64,64)$  & 0.41   & 40   & 4.0795e-4  &  & 0.26   & 29  & 4.3751e-4    \\
$(128,128)$ & 5.42   & 54   & 0.0016     &  & 2.71   & 37  & 0.0018       \\
$(256,256)$ & 27.70  & 73   & 0.0070     &  & 12.53  & 45  & 0.0071       \\
$(512,512)$ & 326.71 & 99   & 0.0288     &  & 135.82 & 49  & 0.0304       \\
 \hline
\end{tabular}
\end{table}
%
\begin{example}\label{Ex3}
Now, we use the nonsymmetric sparse matrix SHERMAN3 of dimension $5005\times 5005$ with $20033$ nonzero entries from the Harwell-Boeing collection \cite{Duff} instead the coefficient matrix $A$. For the coefficient matrix $B$, we use $B={\rm tridiag}(-1,4,-2)$ of dimension $8\times8$ \cite{Khorsand2}.
\end{example}
\begin{table}[h!]
\centering
\caption{\label{Tab3}\small{Results of the Example \ref{Ex3} }}
{\begin{tabular}{@{}lcccc}\hline
Method   &{\centering IT}&{\centering CPU } &{\centering res-norm }  \\
\hline
HSS            & $>10000$ & $>2500$ & 2.32             \\
MRHSS          & $>10000$ & $>2500$ & 1.3021           \\
BiCGSTAB       & $\dag$   & $\dag$  & NaN            \\
HSS-BiCGSTAB   & $\dag$   & $\dag$  & NaN            \\
MRHSS-BiCGSTAB & 12       & 2483.35 & 7.7951e-6     \\
\hline
\end{tabular}}
\end{table}
For this problem, the HSS and the MRHSS methods are converging very slowly. We use the BiCGSTAB method for this problem and observe that this method is diverged. In the Table \ref{Tab3}, dagger $\dag$ shows that no convergence has been obtained. Motivate by \cite{Khorsand2} and \cite{Khorsand3}, we use each of the MRHSS and the HSS methods as a splitting preconditioner in the BiCGSTAB method. We observe that use of the MRHSS method as a precondition improves the results obtained by the corresponding method (MRHSS-BiCGSTAB). However, use of the HSS method as a precondition cannot improve the results.

\section{Conclusion}\label{conc}
In this paper, we have proposed an efficient iterative method, which named the MRHSS method, for solving the continuous Sylvester equation $AX+XB=C$. We have compared the MRHSS method with the HSS method for some problems. We have observed that, for these problems the MRHSS method is more efficient versus the HSS method. Moreover, the use of the MRHSS splitting as a precondition can induce accurate and effective preconditioned BiCGSTAB method.



\end{document}